\documentclass[reqno,11pt]{amsart}
\usepackage{amsmath,amssymb,amsthm}
\usepackage{fullpage}
\usepackage{graphicx}
\usepackage{enumerate}
\usepackage[all]{xy}

\newcommand {\E} {\ensuremath{\mathbb{E}}}

\newcommand {\e} {\ensuremath{\epsilon}}
\newcommand {\g} {\ensuremath{\gamma}}
\newcommand {\G} {\ensuremath{\Gamma}}
\newcommand {\de} {\ensuremath{\delta}}
\newcommand {\lam} {\ensuremath{\lambda}}
\newcommand {\al} {\ensuremath{\alpha}}

\newcommand {\oo} {\ensuremath{\infty}}
\newcommand {\Q} {\ensuremath{\mathbb{Q}}}
\newcommand {\N} {\ensuremath{\mathbb{N}}}
\newcommand {\R} {\ensuremath{\mathbb{R}}}

\newcommand {\Z} {\ensuremath{\mathbb{Z}}}

\newcommand {\bs} {\backslash}

\newcommand{\floor}[1]{\lfloor{#1}\rfloor}
\newcommand{\tr}{\textrm{tr}}
\newcommand{\Stab}{\textrm{Stab}}
\newcommand{\Aut}{\textrm{Aut}}

\newcommand{\spn}{\textrm{span}}
\newcommand {\C} {\ensuremath{\mathbb{C}}}
\newcommand{\scriptP}{\mathcal{P}}

\newtheorem{theorem}{Theorem}
\newtheorem{lemma}{Lemma}[section]
\newtheorem{definition}{Definition}[section]
\newtheorem{cor}{Corollary}[section]
\newtheorem{conjecture}{Conjecture}

\numberwithin{equation}{section}

\title{The rate of escape of random walks on polycyclic and metabelian groups} 
\author{Russ Thompson}
\address{Department of Mathematics, Cornell University, 310 Malott Hall, Ithaca NY 14853, USA} 
\email{thru@math.cornell.edu} 
\date{\today}

\begin{document}

\begin{abstract}
We use subgroup distortion to determine the rate of escape of a simple random walk on a class of polycyclic groups, and we show that the rate of escape is invariant under changes of generating set for these groups. For metabelian groups, we define a stronger form of subgroup distortion, which applies to non-finitely generated subgroups. Under this hypothesis, we compute the rate of escape for certain random walks on metabelian groups via a comparison to the toppling of a dissipative abelian sandpile.	
\end{abstract}

\keywords{law of iterated logarithm, metabelian group, polycyclic group, random walk, rate of escape, abelian sandpile, solvable group, subgroup distortion.}

\subjclass[2000]{60B15, 20F65.}

\maketitle


\section{Introduction}\label{sec:roe_intro}

One of the most basic properties of a random walk is how fast it moves. The first notion one might consider is the \emph{speed} of the random walk,
\begin{align}
\lim_{n\to\oo} \frac{|X_n|_G}{n},
\end{align}
where $|\cdot|_G$ is the word length corresponding to the support of the law of $X_n$.
This limit always exists, and, for non-amenable groups, any simple symmetric random walk has positive speed \cite{Kingman,Derriennic,Guivarch}. A random walk has positive speed if and only if there exist non-constant bounded harmonic functions with respect to the Markov operator associated to the random walk \cite{VaroLongRange}. This also implies that the Poisson boundary of the random walk is non-trivial \cite{KaiVerBoundary}. Thus, for non-amenable groups, the speed is positive. The situation is more complicated and less well understood for amenable groups. Many random walks on amenable groups have zero speed, and so we need different notions to gauge how fast the walk is moving. Below we introduce three such notions, each of them parametrized by an exponent of the number of steps taken. We start with a definition introduced by Revelle \cite{RevelleDrift}.

We recall that a generating set of a finitely generated group is a finite subset of $G$ which is closed under inversion and which generates $G$ as a semigroup. Each generating set of $G$ induces a natural word length on $G$ and these metrics are quasi-isometric. For a particular generating set $S$, we denote this word length by $|\cdot|_S$, but when the choice of generating set does not matter we will use $|\cdot|_G$ to some word length on $G$. We will say a random walk on a finitely generated group $G$ is \emph{adapted} if the law of the random walk, $\mu$, is supported on a finite generating set of $G$ and $\mu(g)=\mu(g^{-1})$ for all $g\in G$. 

\begin{definition}[$\al$-tight degree of escape]
	A random walk has $\al$-tight degree of escape if
	\begin{enumerate}
	\item $\exists \gamma,\delta>0$ such that $P(|X_n|_G>\gamma n^{\al})\geq \delta$,
	\item and $\exists \beta>0$ such that $P(|X_n|_G>xn^{\al})\leq c\exp(-cx^\beta)$ for all $x\geq 0$,
	\end{enumerate}
	hold for all $n\geq 0$.
\end{definition}

The next definition concerns the general behavior of the expected distance travelled by the random walk. For two increasing functions $f,g:\R\to \R$, we write $f\preceq g$ if there exists $c>0$ such that
\begin{align*}
f(x) \leq cg(cx)\ \textrm{ for all $x\in \R$.}
\end{align*}
If $f\preceq g$ and $g \preceq f$ we write $f\simeq g$. We extend this definition to functions from $\N\to\N$ via continuous, piecewise linear extension.

\begin{definition}[Displacement exponent $\al$]
A random walk has displacement exponent $\al$ if 
\begin{align}
\E_\mu|X_n|_G \simeq n^{\al}.	
\end{align}
\end{definition}

It is important to keep the measure $\mu$ driving the random walk in mind when computing this quantity. It is evident that
\begin{align}
	\E_{\mu}|X_n|_S\simeq\E_{\mu}|X_n|_T
\end{align}
for any two distinct word metrics on $G$, but if $\mu$ and $\nu$ are supported on two distinct generating sets of $G$, it is unknown whether or not $\E_{\mu}|X_n|_G$ is coarsely equivalent to $\E_{\nu}|X_n|_G$. This is one of the significant open questions concerning the rate of escape. However, the expected displacement is not a quasi-isometric invariant of graphs due to an example of Benjamini and Revelle \cite{BenjaminiRevelle}.

\begin{definition}[$\al$-law of iterated logarithm]
	A random walk has a law of iterated logarithm with exponent $\al$ if
	\begin{align}
		0<\limsup_{n\to \oo} \frac{|X_n|_G}{ (n\log\log n)^\al } <\oo.
	\end{align}
\end{definition}
In the literature, the above denominator is sometimes called an outer radius for the random walk \cite{GrigoryanEscape} or the upper upper L\'{e}vy class \cite{ReveszBook}.

The methods we use will obtain each of these notions of the rate of escape simultaneously, so we group them under the following definition.

\begin{definition}[Escape exponent $\al$]
An random walk has escape exponent $\al$ if it has $\al$-tight degree of escape, displacement exponent $\al$, and an $\al$-law of iterated logarithm.
\end{definition}

Adapted random walks on groups of polynomial volume growth have escape exponent $1/2$, which follows from the gaussian estimates for the heat kernel \cite{HebischS-C}. For groups of super-polynomial volume growth much less is known. Revelle showed that there exist adapted random walks on Sol and $\textrm{BS}(1,2)$ with escape exponent $1/2$, while Lee and Peres have shown there is a universal lower bound of $1/2$ on the displacement exponent for all adapted random walks on groups \cite{LeePeres}

Examples of groups known to have displacement exponent different than $1/2$ are rare and restricted to iterated lamplighter groups. Independently, Revelle and Erschler showed that there exist adapted random walks on 
\begin{align}
		\Z_2\wr\underbrace{\Z\wr\cdots\wr\Z}_{i\ \textrm{times}}
\end{align}
with displacement exponent $1-2^{-i-1}$ \cite{RevelleDrift,ErschlerDrift}, and Revelle showed that these walks have 
$1-2^{-i-1}$-tight degree of escape \cite{RevelleDrift}. Erschler also showed that iterated wreath products involving $\Z^2$ have expected displacement greater than $n^{1-\e}$ for any $\e>0$ \cite{ErschlerDrift}. This class of examples has been effectively mined out as  $\Z_2\wr\Z^d$, for $d\geq 3$, has non-trivial Poisson boundary for any adapted random walk whose projection to $\Z^d$ is transient \cite{KaiVerBoundary}. Thus these examples have positive speed, or, equivalently, displacement exponent 1.

Knowing the expected displacement of a random walk has applications to other questions concerning the geometry of groups. A lower bound on the expected displacement provides an upper bound on the compression of embeddings of a group into $L^p$ \cite{AustinNaorPeres}, while an upper bound provides a lower bound on Erschler's critical constant for recurrence of $G$-spaces \cite{ErschlerCrt}

We will show that certain polycyclic groups have escape exponent $1/2$ for any adapted random walk. To do this we take advantage of the structure of polycyclic groups in conjunction with the idea of subgroup distortion. 
Subgroup distortion has been much studied in geometric group theory \cite{GerstenDistortion,DavisWreathDistortion}; for polycyclic groups, the primary reference is Osin \cite{OsinExpRad}. 
 
Let $\mathcal{F}$ be a coarse type of functions (the typical examples being polynomial and exponential). We say a subgroup $H < G$ has \emph{upper $\mathcal{F}$ distortion} if there is an invertible function $f$ of coarse type $\mathcal{F}$ such that there exists $c>0$ such that
\begin{align}
|h|_G \leq cf^{-1}(|h|_H)+c,
\end{align}
for all $h\in H$. We will be interested in the case of \emph{upper exponential distortion}, 
\begin{align}\label{defn:SED}
|h|_G \leq c\log(|h|_H+1)+c
\end{align}
for all $h\in H$. Note that upper $\mathcal{F}$ distortion of a group/subgroup pair is does not depend on the choice of word metrics for either the group or its subgroup. For convenience, we assume that the trivial group has every type of upper distortion.

\begin{theorem}\label{thrm:roe_pc}
Let $G$ be a torsion-free polycyclic group satisfying a short exact sequence
\begin{align}\label{eqn:pc_ses}
1\to N\to G\to \Z^r\to1,
\end{align}
where $N$ is finitely generated and nilpotent. If $N$, has upper exponential distortion in $G$, then any simple symmetric random walk on $G$ has escape exponent $1/2$.
\end{theorem}

This theorem applies very broadly to the class of polycyclic groups; by a theorem of Mal'cev \cite{SegalPCBook} every polycyclic group has a torsion-free, polycyclic subgroup of finite index which satisfies \ref{eqn:pc_ses} for some nilpotent group $N$. We will refer to $N$ as the \emph{nilpotent kernel} of $G$.

Upper exponential distortion of the nilpotent kernel is not necessary for a polycyclic group to have escape exponent $1/2$. In particular, the nilpotent kernel of a polycyclic group with polynomial volume growth has upper polynomial distortion, and all simple symmetric random walks on such groups have escape exponent $1/2$. However, the set of elements with upper exponential distortion does not necessarily form a subgroup in polycyclic groups \cite{OsinExpRad}. We examine this phenomenon in abelian-by-cyclic groups in Section \ref{sec:mixed_distortion}.

We also show that certain adapted random walks on metabelian groups have escape exponent $1/2$. Metabelian groups satisfy a short exact sequence
\begin{align}\label{eqn:meta_ses}
1\to A\to G\to \Z^r\to1,
\end{align}
where $A$ is abelian. The abelian kernel of a metabelian group need not be finitely generated, and so upper distortion does not necessarily make sense. We will restrict our attention to cases where the above short exact sequence splits and $r=1$, i.e. $G=A\rtimes_\phi \Z$ for some $\phi\in\Aut(A)$. In this case, the ``distortion" of $A$ may be assessed using properties of $\phi$.

Let $p_\phi$ denote the characteristic polynomial of $\phi$. We define $\max^{(k)}S$ to be the $k$th largest element of a set $S$. We say a polynomial $p\in\Z[t,t^{-1}]$ has property \emph{(EDP)} if 
\begin{align}
	{\max_i}^{(1)} |p_i| \left(1-\frac{1}{ {\max_i}^{(2)}|p_i| } \right)>\sum_{i}|p_i|,	\tag{EDP}
\end{align}
where $p_i$ denotes the coefficient of $t^i$ in $p$.

We will show in Section \ref{sec:sum_met} that if the abelian kernel is finitely generated and $p_\phi$ has a multiple in $\Z[t,t^{-1}]$ which satisfies $(EDP)$, then the abelian kernel has upper exponential distortion. However, when the abelian kernel is not finitely generated, i.e. if $A=\Z[q]^d$, $q\in\Q$, $(EDP)$ can be viewed as a stand in for upper exponential distortion. The first result concerning $(EDP)$ and the rate of escape is the following.

\begin{theorem}\label{thrm:roe_edp}
Let $G=\Z[\rho]^d\rtimes_\phi \Z$ where $\rho$ is an algebraic number. If the characteristic polynomial, $p_\phi$, of $\phi$ has a multiple with property $(EDP)$, then there exist adapted random walks on $G$ with escape exponent $1/2$.
\end{theorem}

We will prove this theorem in Section \ref{sec:topple}. Often, the abelian kernel will not have upper exponential distortion, or, if it is not finitely generated, the characteristic polynomial will not have property $(EDP)$ despite the group having exponential volume growth. This case is treated by the following theorem, which we prove Section \ref{sec:mixed_distortion}.

\begin{theorem}\label{thrm:roe_no_ued}
Let $G=\Z[\rho]^d\rtimes_\phi \Z$ where $\rho$ is an algebraic number. If the characteristic polynomial, $p_\phi$, of $\phi$ can be factored over $\Z$ as $p_+p_0$, where
\begin{enumerate}
\item $p_+$ has a multiple which satisfies $(EDP)$,
\item the roots of $p_+$ are distinct with modulus distinct from 1, and
\item the roots of $p_0$ are distinct with modulus 1,
\end{enumerate}
then there exist adapted random walks on $G$ with escape exponent $1/2$.
\end{theorem}

\textit{Remark.} This result is interesting when $G$ is 2-generated. If $A=A_1\oplus A_2$ where $A_1$ has upper exponential distortion in $G$ and $A_2$ is undistorted in $G$, the rate of escape can be determined for random walks whose law is not supported on any elements with non-trivial projections to both $A_1$ and $A_2$ without appeal to Theorem \ref{thrm:roe_no_ued}.  However, there are cases where the abelian kernel may not split so nicely, but where we can apply Theorem \ref{thrm:roe_no_ued} to determine the rate of escape. Thus we have the following corollary.

\begin{cor}\label{cor:roe_mixed}	
There are 2-generated abelian-by-cyclic groups with exponential volume growth which admit adapted symmetric random walk swith escape exponent $1/2$  and whose abelian kernels have neither upper polynomial nor upper exponential distortion.
\end{cor}

Before delving into the proofs of the above theorems, we present a simple result concerning the relationship between the rate of escape on a group and the rate of escape on quotients of that group. Proofs of the theorems then follow, and we conclude with a discussion of of these results and examples to which they apply.

\section{A comparison theorem for the expected displacement of a random walk}\label{sec:roe_quot}

The principal result of this section is the following lemma. 

\begin{lemma}\label{lem:roe_quotient}
Let $G$ be a finitely generated group with a subgroup $H$, and suppose $X_n$ is a random walk on $G$ driven by a measure $\mu$. Set $\G:=H\bs G$ and let $\pi: G \to \G$ be the canonical projection. Then 
\begin{align}
	\E_\mu|\pi(X_n)|_{\G}\leq \E_\mu|X_n|_G. 
\end{align}
\end{lemma}
\begin{proof}
For each coset $\g \in H\bs G$ we pick a unique $g_\g$ such that $\g=Hg_\g$.  Then we have 
\begin{align}
\E_\mu|\pi(X_n)|_{\G} &= \sum_{\g\in\G}|\g|_\G \mu^{(n)}(\g) \\
				&= \sum_{g_\g: \g\in\G}\min\{|hg_\g|_G\ :\ h\in H\} \mu^{(n)}(Hg_\g) \\
				&= \sum_{g_\g: \g\in\G}\left( \min\{|hg_\g|_G\ :\ h\in H\} \sum_{h\in H} \mu^{(n)}(hg_\g) \right)\\
				&\leq \sum_{g_\g: \g\in\G} \sum_{h\in H}|ng_\g|_G\mu^{(n)}(hg_\g).
\end{align}
For each $g\in G$, there exists $h\in H$ such that $g=hg_\g$ as the cosets of $\G$ partition $G$. Hence, 
\begin{align}
	\sum_{g_\g: \g\in\G} \sum_{h\in H}|hg_\g|_G\mu^{(n)}(hg_\g) &=\sum_{g\in G}|g|_G \mu^{(n)}(g)\\
																&=\E_\mu|X_n|_G. 
\end{align}
\end{proof}

We now consider the potential application of this result to the Hanoi towers groups, $\mathcal{H}$ (see  \cite{NekBook}). The Schreier graph $\G(\mathcal{H},\Stab(1),S)$ is homeomorphic to the infinite Sierpinski graph, and the corresponding limit space for the group is the Sierpinski Gasket. Barlow and Perkins have shown that Brownian motion, $W_t$, on the Sierpinski gasket satisfies
\begin{align}
		\E d(0,W_t)\simeq t^{1/d_w},
\end{align} 
where $d_w= \log(5)/\log(2)$ is the walk dimension \cite{BarlowPerkins}. Teufl has proven a refined version of this fact for a simple symmetric random walk on the infinite Sierpinski graph \cite{Teufl}. Similar estimates may be attainable on many fractal spaces, and one should be able to pass through the homeomorphism between the limit space and corresponding Schreier graph to prove the following conjectures.

\begin{conjecture}
The Hanoi towers group, $\mathcal{H}$, admits a simple symmetric random walk on such that
\begin{align}
		n^{\frac{\log 2}{\log 5}}\preceq \E|X_n|_\mathcal{H}.	
\end{align}
\end{conjecture}

\begin{conjecture}
	Let $G$ be a self-similar group and let $\mathcal{X}$ be its limit space. If $W_t$ is the Brownian motion of $\mathcal{X}$, then for any simple symmetric random walk $X_n$ on $G$,
\begin{align}
		\E d(0,W_n)\preceq \E|X_n|_G.
\end{align}
\end{conjecture}

\section{Upper exponential distortion and the rate of escape}\label{sec:polycyclic}

The goal of this section is to prove the following theorem, from which Theorem \ref{thrm:roe_pc} follows as a corollary when $H$ is nilpotent.

\begin{theorem}\label{thrm:roe_ued}
Let $G$ be a torsion-free, finitely generated group with a short exact sequence,
\begin{align}
	0\to H\to G\to \Z^r\to 0,	
\end{align}
where $H$ is finitely generated and torsion-free. If $H$ has upper exponential distortion in $G$, then any simple symmetric random walk on $G$ has escape exponent $1/2$.
\end{theorem}

The idea underlying this theorem is, that under the assumption of upper exponential distortion of $H$ in $G$, the distance travelled by a random walk after $n$ steps can be controlled by the maximum of that random walk's projection to $\Z^r$. We will prove this theorem by generalizing some results obtained by Pittet and Saloff-Coste \cite{PittetS-Csurvey} to arbitrary generating sets of $G$.

We introduce some notation we will use throughout this section. Fix a finite symmetric generating sets $B$ of $H$ and let $A$ be the canonical basis of $\Z^r$. For a generating set $S$ of $G$, note that the projection $\pi_{\Z^r}(S)$ is a generating set of $\Z^r$, but the projection $\pi_{H}(S)$ is not necessarily a generating set of $H$. This means that $\Z^r$ is undistorted in $G$, but $H$ may be distorted. We will show that $H$ is at most exponentially distorted. 

Since commutators of elements in $\Z^r$ may be non-trivial in $G$, it will be helpful to fix a standard embedding of $\Z^r$ into $G$. Let $k=(k_1,\ldots,k_r)\in\Z^r$, and set $\textbf{k}=a_1^{k_1}\cdots a_r^{k_r}$. The set $K=\{\textbf{k}: k\in\Z^r\}$ is a section of $\Z^r$ in $G$, and every $g\in G$ can be written uniquely in normal form: $g=h\textbf{k}=ha_1^{k_1}\cdots a_r^{k_r}$, $h\in H$, $\textbf{k}\in K$. By $|k|$ we denote the length of $k$ in the canonical basis for $\Z^r$.

We start with an elementary lemma on the expansion factor of automorphisms\cite{PittetS-Csurvey}.
\begin{lemma}\label{lem:aut_estimate}
Let $G$ be a finitely generated group with generating set $S$, and let $H$ be a finitely generated subgroup of $\Aut(G)$ with generating set $T$. Then there exists $q\geq 1$ such that for all $h\in H$ and for all $g\in G$
\begin{align}
		|h\cdot g|_S\leq q^{|h|_T}|g|_S.
\end{align}
\end{lemma}
\begin{proof}
Set $q=\sup_{s\in S, t\in T}|t\cdot s|_S$. The desired estimate follows via induction.
\end{proof}

The following lemma establishes an exponential upper bound on the distortion of $H$ in $G$.

\begin{lemma}\label{lem:length_estimate}
Under the above assumptions, for any generating set $S$ of $G$ there exists non-negative integers $q$ and $C$ such that, given $g=h\textbf{k}$,
\begin{align}
	|h|_B\leq q^{|k|},	
\end{align}
and
\begin{align}
	|k|\leq C|g|_S.	
\end{align}
\end{lemma}

The proof of the lemma relies on the following lemma taken directly from Pittet and Saloff-Coste. We exclude the proof as no changes are necessary to bring it into our present setting. Note that the word metric on $G$ does not appear in this lemma. 

\begin{lemma}\label{lem:comm_estimate}
There exists an integer $q>\max\{|[a,a']|_B\ :\ a,a'\in A\}$ such that for each $\e\in\{-1,1\}$, each $k\in\Z^r$, and each $i\in\{1,\ldots,r\}$, there exists $h\in H$ such that
\begin{align}
	a_1^{k_1}\cdots a_r^{k_r}a_i^\e=ha_1^{k_1}\cdots a_i^{k_i+\e} \cdots a_r^{k_r}
\end{align}
and
\begin{align}
	|h|_B\leq q^{|k|}.	
\end{align}
\end{lemma}

\begin{proof}[Proof of \ref{lem:length_estimate}]
 In the normal form, we can write each $s\in S$ as $s=s_H s_K$, where the factors are the projections to $H$ and $K$. Let $M_H=\max_{s\in H}\{|s_H|_B\}$ and $M_K=\max_{s\in K}\{|s_K|_A\}$. Fix $q$ such that lemmas \ref{lem:comm_estimate} and \ref{lem:aut_estimate} hold and fix $C$ such that $q^C\geq 1+M_H + M_Kq^{M_K}$.
We will induct on $|g|_S$. 

The result is trivial when $|g|_S=0$. We will assume the estimates hold for group element with length at most $l$ relative to $S$. Let $|g|_S=l+1$. Then, for any $s\in S$ such that there exists $g'\in G$ with $g=g's$ and  $|g'|_S=l$, we have by the induction hypothesis
\begin{align}
	g'=h'\textbf{k}'	
\end{align}
where $|h'|_B\leq q^{Cl}$ and $|k'|\leq Cl$. 

Observe that
\begin{align}
	g &= h'\textbf{k}'s_Hs_K \\
	  &= h'(\textbf{k}'s\textbf{k}'^{-1})\textbf{k}'s_K\\	
\end{align}
Set $x=\mathbf{k}'s\mathbf{k}'^{-1}\in H$. By Lemma \ref{lem:aut_estimate}, $|x|_B\leq q^{Cl}|s_H|_B\leq M_H q^{Cl}$. 

We now apply Lemma \ref{lem:comm_estimate} $|s_K|_A$-times to find $y_i\in H$, $i\in\{1,\ldots,|s_K|_A\}$ such that $\textbf{k}'s_K=y_1\cdots y_{|s_K|_A}\textbf{k}$. We have
\begin{align}
	 	|y_i|_B\leq q^{|k'|+i}\leq q^{Cl+M_K},
\end{align}
and $|k|=|k'|+|s_K|_A\leq Cl+M_K\leq C(l+1)$.  

Thus $g=h'xy_1\cdots y_{|s_K|_A}\textbf{k}$, where
\begin{align}
	|h'xy_1\cdots y_{|s_K|_A}|_B &\leq q^{Cl}+M_H q^{Cl}+ M_Kq^{Cl+M_K}\\
								 &\leq q^{C(l+1)},
\end{align}
which completes the proof.
\end{proof}

We next adapt this result to paths in the group.

\begin{lemma}\label{lem:rw_estimate}
Under the conditions of Lemma \ref{lem:length_estimate}, consider a sequence $\sigma=s_1\cdots s_l\in S^l$.
Let
\begin{align}
	k(i)=\pi_{\Z^r}(s_1\cdots s_l)=(k_1(i),\ldots,k_r(i)),
\end{align}
and set
\begin{align}
	M(\sigma)=\max_{1\leq i\leq l}|k(i)|.	
\end{align}
Then there exist constants $q$ and $C$ such that, for any $l$ and any sequence $\sigma\in S^l$,
\begin{align}
	s_1\cdots s_l=h\textbf{k}	
\end{align}
with $h\in H$ and $\textbf{k}\in K$,
\begin{align}
	|h|_B\leq Clq^{C M(\sigma)},	
\end{align}
and $|k|\leq M(\sigma)$.
\end{lemma}
\begin{proof}
Let $M_H=\max_{s\in H}\{|s_H|_B\}$ and $M_K=\max_{s\in K}\{|s_K|_A\}$. Fix $q$ such that Lemmas \ref{lem:comm_estimate} and \ref{lem:aut_estimate} hold and fix $C$ such that $C\geq M_H + M_Kq^{M_K}$.	
	
We will induct on $l$. If $l=0$ the result is trivial. Assume the result holds for any sequence of length $l$. Set $\sigma=(s_1,\ldots,s_{l+1})$ and $\sigma'=(s_1,\ldots,s_{l})$ with $s_i\in S$. Set $g'=s_1\cdots s_l$ and $g=g's_{l+1}$. From the induction hypothesis, 
\begin{enumerate}
\item $g'=h'\textbf{k}'$,
\item $|h'|_B\leq Clq^{C M(\sigma')}$, and
\item $|k'|\leq M(\sigma')$.
\end{enumerate}

Observe that
\begin{align}
	g &= h'\textbf{k}'s_Hs_K \\
	  &= h'(\textbf{k}'s_H\textbf{k}'^{-1})\textbf{k}'s_K.
\end{align}
Set $x=k'sk'^{-1}\in H$. By Lemma \ref{lem:aut_estimate},
\begin{align}
	|x|_B\leq q^{C|k'|}|s_H|_B\leq M_H q^{CM(\sigma')}.
\end{align}

We now apply Lemma \ref{lem:comm_estimate} $|s_K|_A$-times to find $y_i\in H$, $i\in\{1,\ldots,|s_K|_A\}$ such that $\textbf{k}'s_K=y_1\cdots y_{|s_K|_A}\textbf{k}$. We have
\begin{align}
	 	|y_i|_B\leq q^{|k'|+i}\leq q^{CM(\sigma')+M_K}.
\end{align}
It follows from the definition of $M(\sigma)$ and the uniqueness of the normal form that $|k|\leq M(\sigma)$. 

Finally,
\begin{align}
	|h'xy_1\cdots y_{|s_K|_A}|_B &\leq Cl q^{C M(\sigma')}+M_H q^{CM(\sigma')}+ M_Kq^{CM(\sigma')+M_K}\\
								 &= (Cl+M_H+M_Kq^{M_K}) q^{CM(\sigma')}\\
								 &\leq C(l+1)q^{CM(\sigma')}\\
								 &\leq C(l+1)q^{CM(\sigma)}.
\end{align}
\end{proof}

We have now assembled the tools need to prove Theorem \ref{thrm:roe_ued}.

\begin{proof}[Proof of Theorem \ref{thrm:roe_ued}]
	Note that if the path $\sigma$ in Lemma \ref{lem:rw_estimate} is a random walk path in $G$, then $M(\sigma)$ is the maximum of the projection of that random walk to $\Z^r$. We will refer to this random variable as $M_n$. By Lemma \ref{lem:rw_estimate} and the hypothesis we see that if $X_n=h_n\textbf{k}_n$, then
\begin{align}
|\mathbf{k}_n|_S\leq |X_n|_S &\leq |h_n|_S+|\mathbf{k}_n|_S\\
				&\leq C'\log(|h_n|_n+1)+|\mathbf{k}_n|_S\\
				&\leq C(M_n+\log(n))+|\mathbf{k}_n|_S,
\end{align}
for some $C>0$.

To see that $X_n$ has $1/2$-tight degree of escape, first observe that, for a simple symmetric random walk on $\Z^r$, there exist $\gamma,\delta>0$ such that $P(\ |k_n|>\gamma\sqrt(n)\ )\geq\delta$, which establishes the lower bound for $|X_n|_S$. Next, note that there exists $c>0$ such that $P(M_n>x)\leq c'\exp(-c'x^2/n)$ \cite{ReveszBook}, which implies an upper bound of the form
\begin{align}
	P(|X_n|_S>xn^{1/2})\leq c\exp(-cx^2)	
\end{align}
exists for some $c>0$.

The remaining results follow from classical results on simple symmetric random walks and their maxima on $\Z^r$ \cite{ReveszBook}.
\end{proof}

\section{Metabelian groups}\label{sec:roe_meta}

The behavior of the rate of escape for metabelian groups is more complex than that of polycyclic groups. Besides displaying a broader range of known behaviors, metabelian groups are more restrictive in terms of the generating sets for which we can determine the rate of escape. However, by using Lemma \ref{lem:roe_quotient} we can get a rough picture of what rates of escape are possible. This is enabled by the following lemma of Baumslag, which tells us that wreath products essentially serve as universal objects in the class of metabelian groups \cite{BaumslagMetabelian}.

\begin{lemma}\label{lem:baumslag_meta}
Let $G$ be a finitely generated metabelian group. Then there exists a free abelian group $A$ of finite rank and a finitely generated abelian group $T$ such that $G$ is isomorphic to a subgroup of $W/N$, where $W = A\wr T$ and $N$ is a normal subgroup of $W$ contained in $\oplus_T A$.
\end{lemma}

This, along with results from Erschler and Revelle \cite{ErschlerDrift,RevelleDrift}, leads to the following observation.

\begin{lemma}\label{lem:roe_meta_bounds}
Let $G$ be a finitely generated, torsion free metabelian group. Let $W$ be as in Lemma \ref{lem:baumslag_meta}. Then exists an adapted random walk on $G$ such that,
\begin{enumerate}
\item if $d=1$, $n^{1/2}\preceq \E|X_n|\preceq n^{3/4}$,
\item if $d=2$, $n^{1/2}\preceq \E|X_n|\preceq n/\log(n)$, and
\item if $d\geq3$, $n^{1/2}\preceq \E|X_n|\preceq n$,
\end{enumerate}
where $d$ is the dimension of $T$.
\end{lemma}

We will explore the first case of the above lemma more in the following sections. 

\subsection{Abelian-by-cyclic groups and polynomials}\label{sec:abyc_structure}

Let $G=A \rtimes_\phi \Z$, where $A$ is torsion free. We will treat $A$ as a vector space and denote elements thereof in boldface. We will assume that $\phi:\Z\to \textrm{SL}(A)$ is \emph{irreducible}, that is, the characteristic polynomial, $p_\phi$, of $\phi$ is irreducible. This allows us to assume that $G$ is 2-generated by a basis element of $A$ and a generator of $\Z$. Explicitly, these generating sets are of the form
\begin{align}
	S=\{(\pm \mathbf{w},0)\}\cup\{(e,\pm 1)\},
\end{align}
for some basis element $\mathbf{w}$ of $A$. We will concern ourselves only with random walks on such generating sets, but note that one can add elements of the form $(\pm\mathbf{v},0)$ to the generating sets and still apply the techniques described below.

\textit{Remark} If $p_\phi$ were reducible we might need to admit additional basis vectors of $A$ into our generating set. In particular, this happens when $A$ is the direct sum of at least two $\phi$-invariant subspaces. However, it is possible for $\phi$ to be reducible without requiring additional elements in the generating set. This occurs when the direct sum of the $\phi$-invariant subspaces is of finite index in $A$. We treat such an example in Section \ref{sec:sum_mixed_dist}.

Let $X_n=(\textbf{W}_n,Y_n)$ be the random walk on $G$ driven by a symmetric measure $\mu$ on $S$, where $\textbf{W}_n\in A$ and $Y_n\in\Z$. Let $\xi_i=(\mathbf{w}_i,y_i)$ denote the increments of $X_n$, and denote the distributions of $\mathbf{w}_i$ and $y_i$ as $\pi_A(\mu)$ and $\pi_\Z(\mu)$, respectively. Then $Y_n=\sum_{i=1}^{n}y_i$ is a simple symmetric random walk on $\Z$ with distribution $\pi_\Z(\mu)$, but the behavior of $\textbf{W}_n$ is more complicated. Observe that
\begin{align}
\textbf{W}_n	&= \left(\sum_{i=1}^n \phi^{Y_{i-1}}\mathbf{w}_i \right)\\
	&= \left(\sum_{i\in\Z}\omega_i(n)\phi^i\right )\mathbf{w},\label{eqn:poly_form}
\end{align}
where $\omega_i(n)$ are i.i.d. random variables equal in distribution to a simple random walk on $\Z$ with distribution $\pi_\Z(\mu)$, stopped at the random time  
\begin{align}
\theta_i(n)=\#\{0<k\leq n\ |\ Y_k=Y_{k-1}=i \},
\end{align}
which is the local time of $Y_n$. This polynomial can be read as a word representing $\mathbf{W}_n$ in $G$, and $\sum_{i\in\Z}|\omega_i(n)|$ is the length of this word. However, $\sum_{i\in\Z}|\omega_i(n)|=O(n^{3/4})$, and this estimate gives the upper bound for the rate of escape on $\Z\wr\Z$ \cite{RevelleDrift}. When $G$ does not contain a wreath product as a subgroup of finite index this estimate is too large, and we need a way to find a more efficient representation for $\mathbf{W}_n$.

Let $P_n$ denote the random polynomial in $\Z[t,t^{-1}]$ given by \eqref{eqn:poly_form}. The map $\Z[t,t^{-1}] \to A$ corresponds to the composition of the evaluation map $P_n \mapsto P_n(\phi)$ with the action of $\textrm{SL}(A)$ on $A$. 

We now define some properties we will use in the study of $P_n$. The length of a polynomial $p$ is given by $\| p\| _\scriptP:=\sum_i|p_i|$ where $p_i$ is the coefficient of the $i$th degree term of $p$. We let $M(p):=\max\{i\ :\ |p_i|>0\}$ and $m(p):=\min\{i\ :\ |p_i|>0\}$. For convenience we denote $M(P_n)$ by $M_n$ and $m(P_n)$ by $m_n$, and note that these correspond to the maximum and minimum of $Y_n$. We will refer to $d(p)=M(p)-m(p)$ as the \emph{diameter} of $p$. We also set $K(p):=\max_i|p_i|$.  

For $\Z\wr\Z$ each element of $\Z[t,t^{-1}]$ corresponds to a distinct lamp configuration, while in $G$ multiple elements of $\Z[t,t^{-1}]$ may represent the same element of $A$. In particular, if $p_{\phi}$ is the characteristic polynomial of $\phi$, then, by the Cayley-Hamilton theorem, the following diagram commutes.

\begin{displaymath}
    \xymatrix{
        \Z[t,t^{-1}] \ar[dr] \ar[r] & \Z[t,t^{-1}]/\langle p_\phi\rangle \ar[d] \\
             & A }
\end{displaymath}

Hence, rather than using $\| P_n\| _\scriptP$ to estimate $|\textbf{W}_n|_G$, we can obtain more accurate estimates by using a representative in $\Z[t,t^{-1}]/<p_\phi>$ with smaller length. We reduce $P_n$ modulo $p_\phi$ via a process akin to a division algorithm which lowers $K(P_n)$ at the expense of increasing $d(P_n)$.

\subsection{The toppling lemma}\label{sec:topple}

We now restrict our attention solely to the modification of polynomials in $\Z[t,t^{-1}]$. The ability to reduce constant polynomials is also sufficient to establish the upper exponential distortion of the (finitely generated) abelian kernel of an abelian-by-cyclic group, so we present this case separately. Once constant polynomials can be reduced, we can also reduce arbitrary polynomials.

Throughout this section, we will use $\mathit{\log}$ to denote the positive part of the logarithm function, $\mathit{\log}^{(k)}$ to denote the $k$th iterate of the positive part of the logarithm function, and $\mathit{\log}^*$ to denote number of iterates of the logarithm needed to to produce an output less than 1. 

We recall for the reader that a polynomial $p$ has property $(EDP)$ if 
\begin{align}
	{\max_i}^{(1)} |p_i| \left(1-\frac{1}{ {\max_i}^{(2)}|p_i| } \right)>\sum_{i}|p_i|,\tag{EDP}	
\end{align}
where $p_i$ denotes the coefficient of $t^i$ in $p$.

Rather than working with a polynomial $P$ algebraically, we treat the absolute values of the coefficients of $P$ as a height function over $\Z$. This allows us to treat the reduction of $P$ like the toppling of an abelian sandpile (see \cite{RedigMathASM} and the references therein). We can view $\| P \|_\scriptP$ as the total number of grains in the sandpile, and the toppling rule for the sandpile corresponds to adding another polynomial to $P$ when a given coefficient of $P$ is sufficiently large in absolute value. A site in a sandpile topples when it has more grains than a predetermined threshold. For our purposes, we can take the threshold to be some fixed multiple of the largest coefficient of the polynomial we are reducing $P$ by. Note that our toppling rules are spatially homogeneous, i.e. the toppling rules correspond to all $t$-multiples of some polynomial $p$. Thus property $(EDP)$ can be viewed as defining a class of toppling rules.

For the classical abelian sandpile on $\Z$, the toppling rule for site $i$ corresponds to the polynomial $t^{i-1}(t^2-2t+1)$. This toppling rule says that when site $i$ topples it loses two grains and each of its neighbors gains a grain. Note that $t^2-2t+1$ does not have $(EDP)$. This toppling rule corresponds to the Heisenberg group, whose $\Z$-action on $\Z^2$ is given by the matrix $\binom{1\ 1}{0\ 1}$, which has characteristic polynomial $y$. Toppling $P_n$ by $y$ does not produce a sharper estimate on the length of $\mathbf{W}_n$ in this case because $\| P_n \|_\scriptP$ does not change after toppling.

However, if we consider the group Sol instead, whose $\Z$-action on $\Z^2$ is given by the matrix $\binom{2\ 1}{1\ 1}$ which has characteristic polynomial $y(t)=t^2-3t+1$, which has $(EDP)$. In this case, grains of sand are lost during toppling, and so the reduced sandpile will provide a sharper estimate on the length of $\mathbf{W}_n$ in $G$. Such sandpile models where grains of sand can be lost during toppling are called dissipative. 

The following lemmas provide some quantitative bounds on how far a sandpile will spread when it topples according to a polynomial with $(EDP)$.

\begin{lemma}[First toppling lemma]\label{lem:topple}
Fix $K\in \Z$. If $y\in \Z[t,t^{-1}]$ has property $(EDP)$, then there exists $x\in\Z[t,t^{-1}]$ such that $p(t)=K+x(t)y(t)$ satisfies the following
\begin{enumerate}
\item $K(p)<c \max_i|y_i|$ for some $c>0$, and
\item $M(p)-m(p)\leq C(b-a)\log K$ for some $C>0$.
\end{enumerate}
In particular, $\| p\| _\scriptP=O(\log K)$.
\end{lemma}

\begin{proof}
We assume without loss of generality that $K$ is positive and $\max^{(1)}_i|y_i|=y_0>0$. Let $-a$ and $b$ denote the minimal and maximal degrees such that $y_i$ is non-zero. Set $\de:=\| y\| _\scriptP-y_0$ and $r:=\frac{\max_{i\neq0}|y_i|}{y_0}$. Note that $(EDP)$ implies that $\de r<1$.

Set $p^{(1)}=K-\floor{K/y_0}y$. It is clear that $K(p^{(1)})\leq rK$, $M(p^{(1)})=b$, and $m(p^{(1)})=a$. However, the constant term is $K\mod y_0$.

We now proceed by induction on the diameter of $p^{(i)}$, $d(p^{(i)})=M(p^{(i)})-m(p^{(i)})$. We have just done the $d=1$ case. We will keep the salient information about our updates to $p$ in terms of  $K_i:=K(p^{(i)})$ and $d_i:=d(p^{(i)})$.

Suppose we are given $K_i$ and $d_i$. Applying the base case to each term in $p^{(i)}$ yields $d_{i+1}=d_i+a+b$ and
\begin{align}
K_{i+1} &\leq \de r K_i +y_0\\
		&= (\de r)^{i}rK + \sum_{j=0}^{i-1} (\de r)^j y_0\\
		&= (\de r)^{i}rK + \frac{1-(\de r)^i}{1-\de r}y_0.
\end{align}
Thus for $i$ on the order of $\log K$, $K_{i+1}$ will be bounded by some constant multiple of $y_0$, and the desired estimates hold.
\end{proof}

\begin{lemma}[Second toppling lemma]\label{lem:topple_general}
Suppose $y\in \Z[t,t^{-1}]$ has property $(EDP)$. If $P\in\Z[t,t^{-1}]$, then there exists $Q\in\Z[t,t^{-1}]$ such that
\begin{enumerate}
\item $P(A)-Q(A)\in \langle y \rangle$,
\item there exists $C>0$, independent of $P$, such that $K(Q)\leq C$, and
\item $d(Q)\leq d(P) + O\left((\log^*K(P))(\log K(P))\right) $.
\end{enumerate}
\end{lemma}
\begin{proof}
We will assume that $P$ has non-negative coefficients. This is sufficient as, for arbitrary $P$, we have the decomposition $P=Q-R$, where $q_i=p_i$ when $p_i>0$ and zero otherwise and $r_i=p_i$ when $p_i<0$ and zero otherwise. The technique below can then be applied independently to $Q$ and $R$, with the only difference being that in one case we use $-y$ instead of $y$. 	
	
We will apply Lemma \ref{lem:topple} iteratively to $P^{(0)}:=P$. This ensures $(1)$. For the first pass we apply Lemma \ref{lem:topple} to each term of $P$. This produces a polynomial $P^{(1)}$ such that
\begin{enumerate}
\item $K(P^{(1)})\leq |y_0|\log K(P^{(0)})$, and
\item $d(P^{(1)})\leq d(P^{(0)})+C'\log K(P^{(0)})$ for some $C'>0$.
\end{enumerate}
On the $n$th pass we have
\begin{enumerate}
\item $K(P^{(n)})\leq |y_0|\sum_{i=1}^{n-1} \log^{(i)} |y_0|+ \log_n K(P^{(0)})$, and
\item $d(P^{(n)})\leq d(P^{(0)})+C\sum_{i=0}^{n-1}\log K(P^{(i)})$, for some $C>0$.
\end{enumerate}

It takes on the order of $\log^* K(P^{(0)}))$ iterations to bring $K(P^{(n)})$ below some fixed constant $C$ (which depends on $|y_0|$), from which the desired estimates follow.
\end{proof}

We now apply this lemma to prove Theorem \ref{thrm:roe_edp}.

\begin{proof}[Proof of Theorem \ref{thrm:roe_edp}]
Note that
\begin{align}
	|Y_n|_G\leq |X_n|_G\leq |\mathbf{W}_n|_G+|Y_n|_G.	
\end{align}
By dint of the polynomial representation of $\mathbf{W}_n$, Lemma \ref{lem:topple_general} implies that
\begin{align}
	|\mathbf{W}_n|_G &\leq \| Q_n\| _\scriptP \\
						&=O(d(P_n)),
\end{align}
where $Q_n$ is the reduced version of $P_n$. One can observe that $d(P_n)$ is on the order of the maximum of $Y_n$, and so $X_n$ has escape exponent $1/2$ via the reasoning in the proof of Theorem \ref{thrm:roe_ued}.
\end{proof}

\section{Metabelian groups with mixed distortion}\label{sec:mixed_distortion}

We no turn our attention towards proving Theorem \ref{thrm:roe_no_ued}. We will continue to operate under the notation and assumptions introduced in Section \ref{sec:abyc_structure}. We now introduce some notions and notation that will help us talk more accurately about $\Z$-actions on $\Z^d$. 

Fix $\phi\in \textrm{SL}(A)$, and let $d$ denote the abelian rank of $A$. The \emph{expanding eigenspace}, $E_+$, of $\phi$ is the span of the eigenvectors whose eigenvalues have modulus distinct from 1. Note that the expanding eigenspace includes the directions in which $\phi$ is expanding and the directions in which it is contracting. 

The \emph{neutral eigenspace}, $E_0$, of $\phi$ is the span of the eigenvectors whose eigenvalues have modulus 1.  Points in $E_0$ are not necessarily undistorted by the action of $\phi$; for instance, if some of the corresponding eigenvalues have multiplicities greater than 1, we may see polynomial distortion as in the Heisenberg group. Outside of exceptions like this, it is convenient to think of $\phi$ as acting like a rotation on $E_0$. 

We will need to index over the eigenvectors which span these eigenspaces. We will denote their index sets by $I_+$ and $I_0$. We will refer to the eigenvalues as $\{\lam_i\}$ and the corresponding eigenvectors by $\{\mathbf{v}_i\}$. Furthermore, we will assume the eigenvalues satisfy
\begin{align}
	|\lam_1|\geq \cdots \geq |\lam_d|.	
\end{align}

Let $\| \cdot \| $ denote the Euclidean norm on $\C^d$, which projects isometrically to the Euclidean norm on $\R^d$, $\Z^d$, or $\Z[\rho]^d$. Note that the Euclidean metric on $\Z^d$ is equivalent to any intrinsic metric of $\Z^d$ as a finitely generated group. We will also use $\| \cdot\| $ to denote the extension of this norm to linear operators on $\C^d$ (i.e. the $L^2$ norm).

Let $\pi_+$ denote the projection to $E_+$ and $\pi_0$ the projection to $E_0$. We denote the distance of a vector $\mathbf{v}\in \C^d$ from $E_+$ as
\begin{align}
	d_+(\mathbf{v}):=\| (I-\pi_+)\mathbf{v}\| .
\end{align}
Likewise, we denote the distance from $E_0$ as $d_0$.

\begin{lemma}\label{lem:}
There exists $C>0$ such that for all $k\in\Z$ and any $\mathbf{v}\in \C^d$
\begin{align}
	d_0(\phi^k \mathbf{v})\leq C|\lam_1|^k,
\end{align}
and
\begin{align}
	d_+(\phi^k \mathbf{v})\leq C.
\end{align}
\end{lemma}

\begin{proof}
Let $w_i$ denote the coefficients of $\mathbf{w}$ in its eigendecomposition. Then
\begin{align}
	(I-\pi_0)\phi^k\mathbf{w} = \sum_{i\in I_+} \lam_i^k w_i(I-\pi_0) \mathbf{v}_i,
\end{align}
and
\begin{align}
	(I-\pi_+)\phi^k\mathbf{w} = \sum_{i\in I_0} \lam_i^k w_i(I-\pi_+) \mathbf{v}_i.
\end{align}
The desired estimates follow by our assumptions on the moduli of the eigenvalues. 
\end{proof}

The above lemma establishes how far the increments of $\textbf{W}_n$ can move the random walk from $E_0$ and $E_+$,.  This result can be seen as a more precise version of Lemma \ref{lem:aut_estimate} for abelian-by-cyclic groups. When $\Z^d$ has upper exponential distortion in $G$, $E_0=\{0\}$, so in this case the exponential upper bound is sharp.

Using the above lemma, we can see that a random walk will move away from the expanding eigenspace at rate $n^{1/2}$.

\begin{lemma}\label{lem:W}
There exists $C>0$ such that
\begin{align}
	\E_\mu d_+(\textbf{W}_n) \leq C n^{1/2}
\end{align}
for all $n\geq 0$.
\end{lemma}

\begin{proof}
Consider the process $\mathbf{W}^0_n=(I-\pi_+)\textbf{W}_n$, which gives the deviation of $\textbf{W}_n$ from $E_+$. This process is a random walk in $\C^d$ whose increments have bounded length. The increments are of the form 
\begin{align}
	\mathbf{\xi}_{k_n}^\e=\e \phi^{k_n} \mathbf{w},	
\end{align}
where $k_n\in\Z$, $\e\in\{-1,1\}$, and $\mathbf{w}=\sum_{i\in I_0}w_i(I-\pi_0)\mathbf{v}_i$. The distribution on $\e \mathbf{w}$ is uniform, but the distribution of $k_n$ is equal to that of $Y_n$. However, if we let $\mu_n$ denote the distribution on the $n$th increment of the process, we have
\begin{align}
	\sum_{\e}\sum_{k}\mathbf{\xi}_k^\e \mu_n(\xi_k^\e)=0.
\end{align}
As the increments of $\textbf{W}^0_n$ are centered and have bounded length, the desired result holds by an application of the Lindeberg-Feller theorem \cite{Durrett}.
\end{proof}

\begin{figure}
	\begin{center}
\includegraphics[height=3.75in]{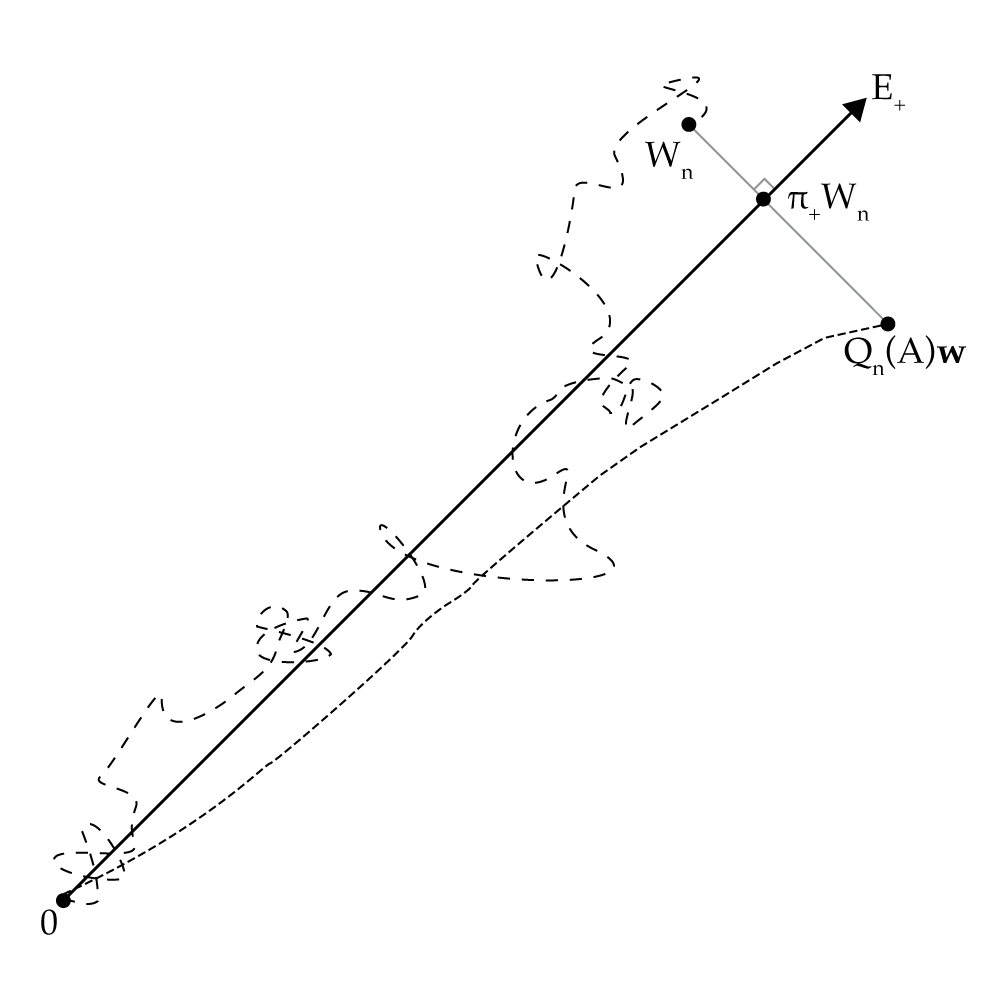}
\caption{The model of an efficient path in $\Z^d \rtimes_A \Z$ projected to $\Z^d$.}
\label{fig:subspace_example}
\end{center}
\end{figure}

We now give a heuristic as to how to take advantage of the above lemma. Suppose $\textbf{W}_n$ is typical for the process. We know it takes on the order of $n^{1/2}$ steps in $A$ to get from $\textbf{W}_n$ to $E_+$. Next we will show that it takes on the order of $n^{1/2}$ steps utilizing the action of $\phi$ to get to a point with the same projection to $E_+$ as $\textbf{W}_n$. In doing so, we will generally pick up some distance from $\mathbf{W}_n$, but so long as this is on the order of $n^{1/2}$ we will have shown that $X_n$ is roughly distance $n^{1/2}$ from the origin. We will do this by using Lemma \ref{lem:topple_general} to reduce $P_n$ relative to $p_{+}:=\sum_{i\in I_+}p_i t^i$.

We will need the following generalization of the Cayley-Hamilton theorem. The proof is a direct computation and is left to the reader.

\begin{lemma}
Let $\phi\in \textrm{SL}_d(\C)$. Denote the eigenvalues of $\phi$ by $\lam_i$ and the corresponding eigenvectors by $v_i$. Fix an index set $I\subseteq\{1,\ldots,d\}$. Let $E_I=\spn_{i\in I}\{v_i\}$. Then $\phi$ is a solution of 
\begin{align}
	p_I(x):=\prod_{i\in I}(x-\lam_i) = 0
\end{align}
over $E_I$.
\end{lemma}

\begin{lemma}\label{lem:Q}
Given $P_n$, if $p_+$ has property $(EDP)$, then there exists $Q_n\in\Z[t,t^{-1}]$ such that 
\begin{enumerate}
\item $P_n-Q_n\in \langle p_{+} \rangle$,
\item $P_n(\phi)-Q_n(\phi)$ acts trivially on $E_+$,
\item $d_+(Q_n(\phi)\mathbf{w})\leq O(\| Q_n\| _\scriptP)$ for all $\mathbf{w}$, and
\item there exists $C>0$, independent of $n$, such that $\| Q_n\| _\scriptP\leq M_n-m_n+C(\log^* K_n )(\log K_n)$.
\end{enumerate}
\end{lemma}

\begin{proof}
We apply Lemma \ref{lem:topple_general} to $P_n$ using $y=p_{+}$. Item $(1)$ follows from the construction used in Lemma \ref{lem:topple_general}, and item $(4)$ is the second conclusion of Lemma \ref{lem:topple_general}.
	
Item $(2)$ follows from $(1)$ and the following observation: if $u\in E_+$, then
\begin{align}
	p_{+}(\phi)\mathbf{u} &=\sum_{i\in I_+}\prod_{i\in I_+} (A-\lam_i I)u_i \mathbf{v}_i\\
				&=0.
\end{align}

For $(3)$ we compute that
\begin{align}
	d_+(Q_n(\phi)\mathbf{w}) &= \| (I-\pi_+)Q_n(\phi)\mathbf{w}\| \\
				&= \| (I-\pi_+)\sum_{i}Q_n(\lam_i)w_i \mathbf{v}_i\| \\
				&= \| \sum_{i\in I_-}Q_n(\lam_i)w_i\mathbf{v}_i\| \\
				&\leq \sum_{i\in I_-}\| Q_n(\lam_i)w_i\mathbf{v}_i\| \\
				&\leq C\sum_{i\in I_-}\| Q_n(\lam_i)\| .
\end{align}
As $\lam_i$ has modulus 1 for each $i\in I_-$, the final sum is on the order of $\| Q_n\| _\scriptP$.
\end{proof}

We can now conclude that the random walk on $G$ has escape exponent $1/2$.

\begin{proof}[Proof of Theorem \ref{thrm:roe_no_ued}]
Recall that $\textbf{W}_n=P_n(\phi)\mathbf{w}$. For the conclusion to hold we only need to show that $|\textbf{W}_n|_G$ has the desired behavior. This is done by showing that $|\textbf{W}_n|_G$ is on the order of $M_n$ plus a term that behaves like the displacement of random walk  with bounded increments on a Euclidean space. As in the proof of Theorem \ref{thrm:roe_pc}, this will be sufficient to imply escape exponent $1/2$. 

Apply Lemma \ref{lem:Q} to obtain $Q_n$. Observe that
\begin{align}
	\textbf{W}_n = Q_n(\phi)\mathbf{w} - (P_n(\phi)-Q_n(\phi))\textbf{w},
\end{align}
and so
\begin{align}
	|\textbf{W}_n|_G &\leq |Q_n(\phi)\mathbf{w}|_G + |(P_n(\phi)-Q_n(\phi))\mathbf{w}|_G.
\end{align}
As $(P_n(\phi)-Q_n(\phi))\mathbf{w}$ lies in $E_0$, we have, by an application of the triangle inequality and Lemmas \ref{lem:W} and \ref{lem:Q},
\begin{align}
	|(P_n(\phi)-Q_n(\phi))\mathbf{w}|_G = O(\| Q_n\| _\scriptP+n^{1/2}).
\end{align}
Combining this with the estimate on $|Q_n(\phi)\mathbf{w}|_G$ to be had from Lemma \ref{lem:Q}, we have
\begin{align}
	|\textbf{W}_n|_G=O(\| Q_n\| _\scriptP+n^{1/2}),
\end{align}
which is on the order of the maximum of $Y_n$. Escape exponent $1/2$ follows from this observation.
\end{proof}

\section{Summary and examples}

\subsection{Polycyclic groups}

Combining Theorem \ref{thrm:roe_pc} with the known escape exponent for groups of polynomial volume growth we have the following.

\begin{theorem}
Let $G$ be a torsion-free polycyclic group. If the nilpotent kernel of $G$ has either upper polynomial distortion or upper exponential distortion then any simple symmetric random walk on $G$ has escape exponent $1/2$. Furthermore, the former case corresponds to $G$ having polynomial volume growth, while the latter corresponds to $G$ having exponential volume growth. 
\end{theorem}

We now provide some tools for determining when the nilpotent kernel of a polycyclic group has upper exponential distortion. Sol and similar groups make for a good place to begin.

\begin{lemma}\label{lem:sed_soltype}
Let $G=\Z^2\rtimes_A \Z$ with $A\in \textrm{SL}_2(\Z)$ with $|\tr(A)|>2$. Then $\Z^2$ is strictly exponentially distorted in $G$.
\end{lemma}
\begin{proof}
Let $B={b_1, b_2}$ be the canonical basis element of $\Z^2$, and let $z$ be the standard basis element of $\Z$. We compute that
\begin{align}
z^k b_i z^{-2k} b_i^{\det(A)^k} &= A^k b_i+\det(A)^k M^-k b_i\\
 				   				&= b_i^{\tr(A^k)}.
\end{align}
The left hand side above has length $3k+2$ in the word metric corresponding to $S$, while $b_i^{\tr(A^k)}$ has length exponential in $k$ in terms on the word metric corresponding to $B$. Thus any $w\in\Z^2$ can be written in $G$ as a word with length on the order of $\log(|w|_B+1)$.		
\end{proof}

We can move beyond abelian-by-cyclic polycyclic groups using the following observation.

\begin{lemma}
Let $G=N\rtimes \Z^r$ be a torsion-free polycyclic group. Consider $N$ as a group of upper-triangular matrices with ones on the diagonal. If two coordinate subgroups, $H_1$ and $H_2$, of $N$ have upper exponential distortion in $G$, then so does $K=[H_1,H_2]$.
\end{lemma}
\begin{proof}
By hypothesis either $H_1\cong H_2 \cong K \cong \Z$ or $K$ is trivial. We have assumed that the trivial group has upper exponential distortion so we will suppose that $K\cong \Z$. Let $k\in K$. Then for any $h_1\in H_1$ and $h_2\in H_2$ such that $k=[h_1,h_2]$,
\begin{align}
|k|_G	&\leq 2(|h_1|_G+|h_2|_G)\\
		&\leq  c_1\log(|h_1|_{H_1}+1) + c_2\log(|h_2|_{H_2}+1)+ c_1 + c_2.
\end{align}
The upper exponential distortion of $K$ follows from this as $|k|_K\geq |h_1|_{H_1}, |h_2|_{H_2}$.
\end{proof}

Combining, the prior lemmas lets us construct nilpotent-by-cyclic groups whose nilpotent kernels have upper exponential distortion.

\begin{cor}
Let $H_3(\Z)$ denote the three dimensional Heisenberg group over $\Z$ and set $G=H_3(\Z)\rtimes_\phi \Z$	where $\phi$ acts on the abelianization of $H_3(\Z)$ as an element of $\textrm{SL}_2(\Z)$, with trace greater than 2 in absolute value. Then $H_3(\Z)$ has upper exponential distortion in $G$. 
\end{cor}

\subsection{Metabelian groups}\label{sec:sum_met}

In what follows we will assume that $G$ is a two-generated abelian-by-cyclic group group. This is done for ease of presentation rather than for any failure of the techniques used. It is clear that we may take any symmetric generating set in the cyclic subgroup without changing the rate of escape as the behavior of the maximum of the random walk on the cyclic subgroup is changed only up to a constant multiple. Admitting additional generators in the abelian kernel is also not problematic. If we have generators $\{\mathbf{w}^{(i)},...,\mathbf{w}^{(k)}\}$, then $W_n=\sum_{i=1}^k P^(i)_n(\phi)\mathbf{w}^{(i)}$, and we now just have to apply the Flattening Lemma independently to each of these polynomials. Thus, so long as our generating set contains no elements with non-trivial projections to both the abelian kernel and the cyclic subgroup, the corresponding random walk will have escape exponent $1/2$. We view this as merely a technical restriction.

The results on metabelian groups apply without hassle to groups whose abelian kernel is $\Z^d$. Next we consider the case when the abelian kernel is of the form $\Z[\rho]^d$, where $\rho$ is an algebraic number. If $\rho$ is rational, then some integer multiple of $p_\phi$ will lie in $\Z[t,t^{-1}]$, so our results readily apply to the Baumslag-Solitar groups $BS(1,n)$ and higher dimensional generalizations thereof. If $\rho$ is not rational, one needs to find a multiple of $p_\phi$ with $(EDP)$ (which implicitly must have integer coefficients).

Since metabelian groups may have non-finitely generated subgroups, subgroup distortion is not the appropriate tool to analyze the behavior of random walks. Instead, for abelian-by-cyclic groups we can generalize to the characteristic polynomial of $\phi$ having property $(EDP)$. 

\begin{lemma}
Let $G=A\rtimes_\phi \Z$ where $A$ is finitely generated. If a multiple of the characteristic polynomial of $\phi$ satisfies $(EDP)$, then $A$ has upper exponential distortion in $G$.
\end{lemma}
\begin{proof}
By the hypothesis, Lemma \ref{lem:topple} states that the polynomial $x(t)=K$ can be represented in $\langle p_\phi \rangle$ by a polynomial $p$ with $\| p\| _\scriptP=O(\log K)$. Via the polynomial model for $G$, $p$ represents a word in $G$ with length $\| p\| _\scriptP$. This implies that any $d$-tuple in $A$, where $d$ is the abelian rank of $A$, can be represented in $G$ by a word whose length is logarithmic in terms of the length in $A$, and so upper exponential distortion holds.
\end{proof}

We now highlight a specific case of the above lemma. For $\phi\in \textrm{SL}_2(\Z)$, there are hyperbolic $\phi$ with $|\tr(\phi)|\leq 2$. One can check that in this case $\det(\phi)=-1$. Then
\begin{align}
		(x^2+\tr(\phi)x-1)p_\phi(x)=x^4 -(2+\tr(\phi)^2)x^2+1,
\end{align} 
so these examples have upper exponential distortion as well. A similar analysis to this was done by Warshall to study the existence of dead-ends in the Cayley graphs of abelian-by-cyclic groups \cite{WarshallDeepPockets}. This observation, along with the fact that non-hyperbolic automorphisms imply polynomial volume growth, leads to the following corollary

\begin{cor}\label{cor:2d_ued}
Let $\phi$ be a hyperbolic automorphism of $\Z^2$. Then $\Z^2$ has upper exponential distortion in $\Z^2\rtimes_\phi \Z$ if and only if $\phi$ is hyperbolic.
\end{cor}

For groups with higher dimensional abelian kernels the situation is more complicated. We treat this case in the following section. We close this section with an example where $G$ is not 2-generated.

\emph{Example.} Let $\lam=\frac{1+\sqrt{5}}{2}$, and note that $1/\rho=\rho-1$. Consider the group $G=\Z[\rho]^2\rtimes_\phi\Z$, where 
\begin{align}
	\phi=\begin{pmatrix}
	\lam & 0\\
	0 & \lam^{-1}\\
	\end{pmatrix}.
\end{align}
Since $\phi$ is diagonal we need to take two generators from $\Z[\rho]^2$, namely $\mathbf{e}_1$ and $\mathbf{e}_2$, the standard basis vectors of $\Z^2$. Note that $\Z[\rho]^2=\spn\{\mathbf{e}_1,\mathbf{e}_2,\sqrt{5}\mathbf{e}_1,\sqrt{5}\mathbf{e}_2\}$. The characteristic polynomial of $\phi$ is $x^2-\sqrt{5}x+1$, which does not lie in $\Z[t,t^{-1}]$. However,
\begin{align}
	(x^2+\sqrt{5}x+1)p_\phi=x^4-3x^2+1,	
\end{align}  
so property $(EDP)$ is satisfied, and we can conclude that any symmetric random walk on $G$ whose law has support $\{(\pm \mathbf{e}_1,0),(\pm \mathbf{e}_2,0),(\mathbf{0},\pm 1)\}$ has escape exponent $1/2$.

It is interesting to note that $\Z[\rho]^2$ is isomorphic to $\Z^4$, and that the action induced on $\Z^4$ by $\phi$ corresponds to the matrix (with coordinates ordered as $(\mathbf{e}_1,\sqrt{5}\mathbf{e}_1,\mathbf{e}_2\sqrt{5}\mathbf{e}_2)$),
\begin{align}
M=\begin{pmatrix}
0 & 1 & 0 &0\\
1 & 1 & 0 & 0\\
0 & 0 & -1 & 1\\
0 & 0 & 1 & 0\\
\end{pmatrix},
\end{align}
which has characteristic polynomial $p_M (x)=x^4-3x^2+1$. The moral of this example is that if one wishes to consider abelian kernels of the form $\Z[\rho]^d$ where $\rho$ is an algebraic number, it may pay rewrite the group in terms of an action on $\Z^{d'}$ for suitably chosen $d'$.

\subsection{Abelian kernels without upper exponential distortion}\label{sec:sum_mixed_dist}

There are two examples we wish to highlight. The first example is due to Conner \cite{Conner}. Let $G_1=\Z^4\rtimes_\phi\Z$ where
\begin{align}
\phi=\begin{pmatrix}
2 & -1 & 2 &-1\\
1 & 0 & 0 & 0\\
0 & 1 & 0 & 0\\
0 & 0 & 1 & 0\\
\end{pmatrix}.
\end{align}
$G$ has exponential volume growth, but every cyclic subgroup is undistorted in $G$ \cite{Conner, OsinExpRad}.

The matrix $\phi$ has two positive, real eigenvalues distinct from 1 and two complex eigenvalues of modulus 1.  The characteristic polynomial factors as
\begin{align}
	p_\phi(x)=(x^2-(1+\sqrt{2})x+1)(x^2-(1-\sqrt{2})x+1),	
\end{align}
where the first factor corresponds to the eigenvalues with roots of modulus distinct from 1. That the characteristic polynomial does not factor over $\Z$ makes this example difficult to work with.

We also consider $G_2=\Z^4\rtimes_\psi \Z$ for 
\begin{align}
	\psi=\begin{pmatrix}
	3 & -2 & 3 & -1\\
	1 & 0 & 0 & 0\\
	0 & 1 & 0 & 0\\
	0 & 0 & 1 & 0\\
	\end{pmatrix}.
\end{align}	
This group has two real eigenvalues distinct from $1$ and two complex eigenvalues of modulus 1, namely $i$ and $-i$. The characteristic polynomial factors as 
\begin{align}
		p_\psi(x)=(x^2-3x+1)(x^2+1),
\end{align}
The first factor corresponds to the eigenvalues with modulus distinct from 1 while the other factor corresponds to the eigenvalues with modulus 1.

Every cyclic subgroup of $G_1$ is undistorted \cite{Conner}. This is because the $E_+$ and $E_0$ for $\phi$ only intersect $\Z^4$ at the origin. However, the characteristic polynomial of $\psi$ splits over $\Z$ so $E_+$ and $E_0$ have non-trivial intersections with $\Z^4$. In particular, $(E_+ \oplus E_0)\cap \Z^4$ has finite index in $\Z^4$. One can also check that $G_2$ is two generated if $\mathbf{w}$ is taken to be one of the canonical basis vectors of $\Z^4$; this happens because the basis vectors do not lie in $(E_+ \oplus E_0)\cap \Z^4$. This establishes Corollary \ref{cor:roe_mixed}.

\bibliography{Research}
\bibliographystyle{alpha}

\end{document}